\input amstex
\documentstyle{amsppt}
\magnification=\magstep1                        %<====
\hsize6.5truein\vsize8.9truein                  %<====
\NoRunningHeads
\loadeusm

\magnification=\magstep1                        %<====
\hsize6.5truein\vsize8.9truein                  %<====
\NoRunningHeads
\loadeusm

\document
\topmatter

\title
Reverse Markov- and Bernstein-type inequalities for incomplete polynomials
\endtitle

\rightheadtext{reverse Markov- and Bernstein-type inequalities}

\author Tam\'as Erd\'elyi
\endauthor

\address Department of Mathematics, Texas A\&M University,
College Station, Texas 77843, College Station, Texas 77843 \endaddress
\endthanks


\dedicatory Dedicated to Professor Richard Varga on the occasion of his 90th birthday
\enddedicatory

\date September 14, 2018
\enddate

\email terdelyi\@math.tamu.edu
\endemail

\abstract
Let ${\Cal P}_k$ denote the set of all algebraic polynomials of degree at most $k$ with real coefficients.
Let ${\Cal P}_{n,k}$ be the set of all algebraic polynomials of degree at most $n+k$ having exactly $n+1$ zeros at $0$.  
Let 
$$\|f\|_A := \sup_{x \in A}{|f(x)|}$$ 
for real-valued functions $f$ defined on a set $A \subset {\Bbb R}$. Let 
$$V_a^b(f) := \int_a^b{|f^{\prime}(x)| \, dx}$$
denote the total variation of a continuously differentiable function $f$ on an interval $[a,b]$.  
We prove that there are absolute constants $c_1 > 0$ and $c_2 > 0$ such that
$$c_1 \frac nk\leq \min_{P \in {\Cal P}_{n,k}}{\frac{\|P^{\prime}\|_{[0,1]}}{V_0^1(P)}}  
\leq \min_{P \in {\Cal P}_{n,k}}{\frac{\|P^{\prime}\|_{[0,1]}}{|P(1)|}} \leq c_2 \left( \frac nk + 1 \right)$$ 
for all integers $n \geq 1$ and $k \geq 1$. We also prove that there are absolute constants $c_1 > 0$ and $c_2 > 0$ such that 
$$c_1 \left(\frac nk\right)^{1/2} \leq \min_{P \in {\Cal P}_{n,k}}{\frac{\|P^{\prime}(x)\sqrt{1-x^2}\|_{[0,1]}}{V_0^1(P)}} \leq 
\min_{P \in {\Cal P}_{n,k}}{\frac{\|P^{\prime}(x)\sqrt{1-x^2}\|_{[0,1]}}{|P(1)|}} \leq c_2 \left(\frac nk + 1\right)^{1/2}$$ 
for all integers $n \geq 1$ and $k \geq 1$.

\endabstract

\endtopmatter

\head 1. Introduction and Notation \endhead

In April, 2018, A. Eskenazis and P. Ivanisvili [8] asked me if I knew polynomial inequalities of a certain type. 
The inequalities they were interested in looked to me immediately as reverse (or inverse) 
Markov- and Bernstein-type inequalities for incomplete polynomials on the interval $[0,1]$, 
but I have not been aware of any such inequalities in the literature. This short paper 
is a result of an effort to answer the questions raised by A. Eskenazis and P. Ivanisvili [8].   
G.G. Lorentz, M. von Golitschek, and Y. Makovoz devotes Chapter 3 of their book [12] to incomplete polynomials. 
E.B. Saff and R.S. Varga were among the researches having contributed significantly to this topic. 
See [17] and [18], for instance. See also [1] written by I. Borosh, C.K. Chui, and P.W. Smith.  
Reverse Markov- and Bernstein type inequalities were first studied by P. Tur\'an [19] and J. Er\H od [7] 
in 1939. The research on Tur\'an and Er\H od type reverses of Markov- and Bernstein-type inequalities 
suddenly got a new impulse in 2006 in large part by the work of Sz. R\'ev\'esz [15], and 
several paper have been publishes on such inequalities in recent years, see [5], [6], [9], [10], [11], [13], [14], 
[16], [20], and [21], for example.         

Let ${\Cal P}_k$ denote the set of all algebraic polynomials of degree at most $k$ with real coefficients.
Let ${\Cal P}_{n,k}$ be the set of all algebraic polynomials of degree at most $n+k$ having exactly $n+1$ zeros at $0$.
That is, every $P \in {\Cal P}_{n,k}$ is of the form
$$P(x) = x^{n+1}R(x)\,, \qquad R \in {\Cal P}_{k-1}\,.$$ 
Let 
$$\|f\|_A := \sup_{x \in A}{|f(x)|}$$ 
for real-valued functions $f$ defined on a set $A \subset {\Bbb R}$. Let
$$V_a^b(f) := \int_a^b{|f^{\prime}(x)| \, dx}$$
denote the total variation of a continuously differentiable function $f$ on an interval $[a,b]$. 

\head 2. New Results \endhead

\proclaim{Theorem 2.1} There are absolute constants $c_1 > 0$ and $c_2 > 0$ such that
$$c_1 \frac nk \leq \min_{P \in {\Cal P}_{n,k}}{\frac{\|P^{\prime}\|_{[0,1]}}{V_0^1(P)}}  
\leq \min_{P \in {\Cal P}_{n,k}}{\frac{\|P^{\prime}\|_{[0,1]}}{|P(1)|}} \leq c_2 \left( \frac nk +1 \right)$$
for all integers $n \geq 1$ and $k \geq 1$. Here $c_1 = 1/12$ is a suitable choice.
\endproclaim

\proclaim{Theorem 2.2} There are absolute constants $c_1 > 0$ and $c_2 > 0$ such that
$$c_1 \left(\frac nk\right)^{1/2} \leq \min_{P \in {\Cal P}_{n,k}}{\frac{\|P^{\prime}(x)\sqrt{1-x^2}\|_{[0,1]}}{V_0^1(P)}} \leq 
\min_{P \in {\Cal P}_{n,k}}{\frac{\|P^{\prime}(x)\sqrt{1-x^2}\|_{[0,1]}}{|P(1)|}} \leq c_2 \left(\frac nk + 1\right)^{1/2}$$
for all integers $n \geq 1$ and $k \geq 1$. Here $c_1 = 1/6$ is a suitable choice. 
\endproclaim

\head 3. Lemmas \endhead

Our first lemma is a simple consequence of the well known Chebyshev's inequality (see p. 235 of [4], for instance) 
on the growth of polynomials.

\proclaim {Lemma 3.1}
We have
$$|Q(x)| \leq |2x|^k \|Q\|_{[-1,1]}, \qquad x \in {\Bbb R} \setminus (-1,1)\,,$$
for every $Q \in {\Cal P_k}$, $k \geq 0$. 
\endproclaim

The following lemma follows from Lemma 3.1 by a simple linear transformation.

\proclaim {Lemma 3.2}
Let $a,b \in {\Bbb R}$ and $a < b$. We have
$$|Q(x)| \leq \left|\frac{4x-2(a+b)}{b-a}\right|^{k} \|Q\|_{[a,b]}\,, \qquad x \in {\Bbb R} \setminus (a,b)\,,$$
for every $Q \in {\Cal P_k}$, $k \geq 0$.
\endproclaim

Our next lemma is a special case of Lemma 3.2.

\proclaim {Lemma 3.3}
Let $a,b \in {\Bbb R}$ and $a < b$. We have
$$|Q(x)^2(1-x^2)| \leq \left|\frac{4x-2(a+b)}{b-a}\right|^{2k} \|Q(u)^2(1-u^2)\|_{[a,b]}, \qquad x \in {\Bbb R} \setminus (a,b)\,,$$
for every $Q \in {\Cal P_{k-1}}$, $k \geq 1$. 
\endproclaim

We will use Lemmas 3.2 and 3.3 to prove our next couple of lemmas.

\proclaim {Lemma 3.4}
Let $n \geq 1$ and $k \geq 1$ be integers, and let $S(x) := x^nR(x)$ with $R \in {\Cal P}_k$.
We have
$$|S(x)| \leq x^{n/2} \|S\|_{[0,1]}\,, \qquad x \in [0,1-10k/n]\,.$$
\endproclaim

\proclaim {Lemma 3.5}
Let $n \geq 1$ and $k \geq 1$ be integers, and let $S(x) := x^nQ(x)\sqrt{1-x^2}$ with $Q \in {\Cal P}_{k-1}$.
We have  
$$|S(x)| \leq x^{n/2} \|S\|_{[0,1]}\,, \qquad x \in [0,1-10k/n]\,.$$
\endproclaim

To prove the upper bounds in Theorems 2.1 and 2.2 we need the following result proved in [2].

\proclaim{Lemma 3.6}
Let $\nu \geq 0$ and $\kappa \geq 1$ be nonnegative integers. There is an absolute constant $c_3 > 0$ such that
$$|P^{\prime}(x)| \leq c_3 \left( \frac{(\nu+\kappa)\kappa}{x(1-x)} \right)^{1/2} \|P\|_{[0,1]}\,, \qquad x \in (0,1)\,,$$
for every polynomial $P \in {\Cal P}_{\nu+\kappa}$ having at most $\kappa$ zeros in the open disk with diameter $(0,1)$.
\endproclaim

Let $k \geq 1$ be an integer. We define
$$\alpha_j = 1 + \cos \left(\pi \frac{2j-1}{4k} \right)\,, \qquad j=1,2,\ldots,k\,.$$
Let $n = 2k + m$, where $m \geq 1$ is an integer. Let $1 > \gamma_1 > \gamma_2 >\cdots > \gamma_k > 0$ be defined by
$$\gamma_j := \frac{\alpha_j - (1-m/k)}{1+m/k}\,, \qquad j=1,2,\ldots,k\,.$$ 
Let $q_n \in {\Cal P_n}$ be the unique polynomial of the form 
$$q_n(x) = (x+1)^{n-k}\prod_{j=1}^k{(x-\rho_j)}$$
equioscillating $k+1$ times on $[-1,1]$ between $-1$ and $1$, that is, there are 
$$1 = x_0 > x_1 > \cdots > x_k > -1$$  
satisfying
$$q_n(x_j) = (-1)^j = (-1)^j\|q_n\|_{[-1,1]}\,, \qquad j=0,1,\ldots,k\,.$$

To prove the upper bounds in Theorems 2.1 and 2.2 we also need the following lemma 
stating a key observation from the proof of Lemma 4 in [2]. 

\proclaim{Lemma 3.7}
With notation introduced above we have $\rho_j \leq \gamma_j$ for each $j=1,2,\ldots,k$. 
As a consequence, there is an absolute constant  $c_4 > 0$ such that
$$\rho_j \leq 1 - \frac{c_4j^2}{nk}\,, \qquad j=1,2,\ldots,k\,.$$
\endproclaim

For our purpose to prove the upper bounds in Theorems 2.1 and 2.2 the following version 
of Lemma 3.7 will be convenient for us.

\proclaim{Lemma 3.8}
Let $1 \leq \kappa \leq \nu-1$ be integers. Let $T := T_{\nu,\kappa}$ be the Chebyshev polynomial for the
M\"untz space
$$\text{\rm span}\{x^{\nu}, x^{\nu+1}, \ldots, x^{\nu+\kappa}\}$$
on $[0,1]$ normalized so that $T(1)=1$. Denote the zeros of $T$ in $(0,1)$ by 
$$\beta_1 > \beta_2 > \cdots > \beta_{\kappa}\,.$$
We have
$$\beta_j \leq 1 - \frac{c_4j^2}{(\nu+\kappa)\kappa} 
\leq 1 - \frac{c_4j^2}{\nu\kappa} \,, \qquad j=1,2,\ldots,\kappa\,,$$
where $c_4 > 0$ is the absolute constant appearing in Lemma 3.7.
\endproclaim

In fact, what we need in the proofs of the upper bounds in Theorems 2.1 and 2.2 is the following easy 
consequence of Lemma 3.8.

\proclaim{Lemma 3.9}
Let $\kappa \geq 2$ and $20\kappa \leq \nu$ be integers. Let $T := T_{\nu,\kappa}$ be the Chebyshev polynomial for the M\"untz space 
$$\text{\rm span}\{x^{\nu}, x^{\nu+1}, \ldots, x^{\nu+\kappa}\}$$
on $[0,1]$ normalized so that $T(1)=1$. There is an absolute constant $c_5 > 0$ such that
$$\int_0^1{T(u)^2 \, du} \geq \frac{c_5\kappa}{\nu}\,.$$
\endproclaim

\head 4. Proofs of the Lemmas \endhead

\demo {Proof of Lemma 3.1}
Let $T_k$ be the $k$-th Chebyshev polynomial defined by
$$T_k(\cos t) = \cos(kt)\,, \qquad t \in {\Bbb R}\,.$$
It is well known that
$$T_k(x) = 2^{k-1} \prod_{j=1}^k{(x-x_j)}$$
where
$$x_j = \cos\left(\frac{(2j-1)\pi}{2k}\right)\,, \qquad j=1,2,\ldots,k\,,$$
and hence
$$1 > x_1 > x_2 > \cdots > x_k > -1\,.$$
Using Chebyshev's inequality (see E.2 on page 235 of [4], for instance) we have
$$\split |Q(x)| \leq & |T_k(x)| \dot \|Q\|_{[-1,1]} = \left( 2^{k-1} \prod_{j=1}^k{|x-x_j|} \right) \|Q\|_{[-1,1]} \cr 
= & \left( 2^{k-1} \prod_{j=1}^k{|x^2-x_j^2|^{1/2}} \right) \|Q\|_{[-1,1]} \leq |2x|^k \|Q\|_{[-1,1]} \cr \endsplit$$
for every $Q \in {\Cal P_k}$ and $x \in {\Bbb R} \setminus (-1,1)$.
\qed \enddemo

\demo {Proof of Lemma 3.4}
Let $\delta := k/n \in (0,1)$. Without loss of generality we may assume that $\|S\|_{[0,1]} = 1$. Therefore
$$\|R\|_{[1-\delta,1]} \leq (1-\delta)^{-n}\,.$$
Combining this with Lemma 3.3 we obtain that if $x \in [0,1-\delta]$, then 
$$\split |S(x)| \leq & \, x^{n/2} \cdot |x^{n/2}R(x)| \cr
\leq & \, x^{n/2} \cdot x^{n/2} \cdot \left|\frac{4x-(4-2\delta)}{\delta}\right|^k \|R\|_{[1-\delta,1]} \cr  
\leq & \, x^{n/2} \cdot x^{n/2} \cdot \left(\frac{4-4x}{\delta}\right)^k (1-\delta)^{-n} \cr 
\leq & \, x^{n/2} \cdot x^{n/2}\left(\frac{4-4x}{\delta}\right)^k (1-\delta)^{-n} = x^{n/2} f(x)\,, \cr \endsplit$$
where
$$f(x) = x^{n/2}\left(\frac{4-4x}{\delta}\right)^k (1-\delta)^{-n}\,. \tag 4.1$$ 
To finish the proof we need to show that $|f(x)| \leq 1$ for every $x \in [0,1-10k/n]$. 
If $1-10k/n < 0$ there is nothing to prove, so we may assume that $1-10k/n \geq 0$.  
The function is clearly nonnegative on $[0,1]$, and by examining the sign of 
$f^{\prime}(x)$ it is easy to see that $f$ is increasing on the interval $[0,n/(n+2k)]$,  
and hence on $[0,1-2k/n] \subset [0,n/(n+2k)]$ as well. 
Using (4.1) to estimate the value of $f$ at $x = 1-10k/n \geq 0$,  
we obtain
$$\split f(x) = & (1-10k/n)^{n/2} 40^k (1-k/n)^{-n} \leq (1-5k/n)^n 40^k (1-k/n)^{-n} \cr 
\leq & (1-4k/n)^n 40^k \leq e^{-4k}40^k = \left( \frac{40}{e^4} \right)^k \leq 1\,, \cr \endsplit$$
hence $0 \leq f(x) \leq 1$ for every $x \in [0,1-10k/n]$, indeed.
\qed \enddemo   

\demo {Proof of Lemma 3.5}
Applying Lemma 3.4 with $S \in {\Cal P}_{2n}$ defined by $S(x)^2 = x^{2n}U(x)$, where $U \in {\Cal P}_{2k}$ is defined 
by $U(x) = Q(x)^2(1-x^2)$, we obtain the lemma. 
\qed \enddemo

\demo{Proof of Lemma 3.9}
As before, denote the zeros of $T$ in $(0,1)$ by 
$$\beta_1 > \beta_2 > \cdots > \beta_{\kappa}\,.$$
We introduce the points of equioscillation $x_0 > x_1 > \cdots > x_{\kappa}$, that is,
$T(x_j) = (-1)^j$ and $\beta_j \in (x_j,x_{j-1})$ for $j=1,2,\ldots,\kappa$,
where $x_{\kappa} \geq 1 - 10\kappa/\nu \geq 1/2$ follows from Lemma 3.4 and the assumption $20\kappa \leq \nu$.
We define  $y_j \in (\beta_{j+1},x_j)$ by
$$T(y_j) = (-1)^j(1/2)\,, \qquad j=1,2,\ldots,\kappa-1\,.$$
The Mean Value Theorem and Lemma 3.6 imply that there are a $\xi_j \in (x_j,y_j)$ such that
$$\split 1/2 = & |T(x_j) - T(y_j)| = (x_j - y_j)|T^{\prime}(\xi_j)| 
\leq (x_j - y_j) c_3 \left( \frac{(\nu+\kappa)\kappa}{\xi_j(1-\xi_j)} \right)^{1/2} \cr 
\leq & c_3(x_j - y_j) \left( \frac{(\nu+\kappa)\kappa}{(1/2)(1-\beta_j)} \right)^{1/2}\,, \qquad j=1,2,\ldots,\kappa-1\,, \cr \endsplit$$
and hence
$$x_j - y_j \geq c_6 \frac{(1-\beta_j)^{1/2}}{(\nu\kappa)^{1/2}}\,, \qquad j=1,2,\ldots,\kappa-1\,,$$
with an absolute constant $c_6 > 0$. 
Observe that $|T(x)| \geq 1/2$ on each of the intervals $[y_j,x_j]$, $j=1,2,\ldots,\kappa-1$, so
$$m(\{x \in [0,1]:|T(x)| \geq 1/2\}) \geq \sum_{j=1}^{\kappa-1}{(x_j - y_j)} 
\geq \sum_{j=1}^{\kappa-1}{c_6\frac{(1-\beta_j)^{1/2}}{(\nu\kappa)^{1/2}}}\,,$$
where $m(A)$ denotes the Lebesgue measure of a set $A \subset {\Bbb R}$. Combining this with
Lemma 3.8, we obtain
$$m(\{x \in [0,1]:|T(x)| \geq 1/2\}) \geq \sum_{j=1}^{\kappa-1}{c_6\frac{(c_4 j^2/(\nu\kappa))^{1/2}}{(\nu\kappa)^{1/2}}} 
\geq c_7 \sum_{j=1}^{\kappa-1}{\frac{j}{\nu\kappa}} = \frac{c_5\kappa}{\nu}$$ 
with some absolute constants $c_7 > 0$ and $c_5 > 0$, and the lemma follows.
\qed \enddemo

\head 5. Proofs of Theorems 2.1 and 2.2.\endhead

\demo {Proof of the lower bound in Theorem 2.1}
Let $P \in {\Cal P}_{n,k}$ be of the form
$$P(x) = x^{n+1}R(x)\,, \qquad R \in {\Cal P}_{k-1}\,.$$
We define
$$S(x) := P^{\prime}(x) = x^nQ(x)\,,$$
where $Q \in {\Cal P}_{k-1}$ is defined by
$$Q(x) = (n+1)R(x) - xR^{\prime}(x)\,.$$
We also define $y := \max\{1-10k/n,0\} \geq 0$.
We have
$$V_0^1(P) = \int_0^1{|P^{\prime}(x)| \, dx} = 
\int_0^y{|P^{\prime}(x)| \, dx} + \int_y^1{|P^{\prime}(x)| \, dx} \tag 5.1$$
The first term at the right-hand side of (5.1) can be estimated by Lemma 3.4 as
$$\int_0^y{|P^{\prime}(x)| \, dx} = \int_0^y{|S(x)|\, dx} 
\leq \int_0^y{(x^{n/2} \|S\|_{[0,1]})\, dx}  \leq \frac 2n \,\|S\|_{0,1]}i\,, \tag 5.2$$
while the second term at the right-hand side of (5.1) can be estimated as 
$$\int_y^1{|P^{\prime}(x)| \, dx} = \int_y^1{|S(x)|\, dx} \leq (1-y) \|S\|_{[0,1]} 
\leq \frac{10k}{n} \, \|S\|_{[0,1]}\,.  \tag 5.3$$
Combining (5.1), (5.2), and (5.3) we obtain
$$V_0^1(P) =  \int_0^1{|P^{\prime}(x)| \, dx} \leq \frac 2n \|S\|_{[0,1]} + \frac{10k}{n} \, \|S\|_{[0,1]}  
\leq \frac{10k+2}{n} \, \|P^{\prime}\|_{[0,1]}\,,$$
and the lower bound of Theorem 2.1 follows.
\qed \enddemo

\demo {Proof of the lower bound in Theorem 2.2}
We define 
$$S(x) := P^{\prime}(x)\sqrt{1-x^2} = x^nQ(x)\sqrt{1-x^2}\,,$$
where $Q \in {\Cal P}_{k-1}$ is defined by  
$$Q(x) = (n+1)R(x) - xR^{\prime}(x)\,.$$
We also define $y := \max\{1-10k/n,0\} \geq 0$, as in the proofs of Theorems 2.1. 
We have
$$V_0^1(P) = \int_0^1{|P^{\prime}(x)| \, dx} = 
\int_0^{y}{|P^{\prime}(x)| \, dx} + \int_{y}^1{|P^{\prime}(x)| \, dx} \,. \tag 5.4$$
The first term at the right-hand side of (5.4) can be estimated by Lemma 3.5 as
$$\split \int_0^y{|P^{\prime}(x)| \, dx} = & \int_0^y{|S(x)|(1-x^2)^{-1/2} \, dx} \cr
\leq & \int_0^y{(x^{n/2} \|S\|_{[0,1]})(1-y^2)^{-1/2} \, dx} \cr 
\leq & \frac 2n \, \|S\|_{[0,1]}(1-y)^{-1/2} \leq \frac 2n \, \|S\|_{[0,1]}(10k/n)^{-1/2} \cr 
\leq & (kn)^{-1/2} \|S\|_{0,1]}\,, \cr \endsplit \tag 5.5$$
while the second term at the right-hand side of (5.4) can be estimated as
$$\split \int_y^1{|P^{\prime}(x)| \, dx} = & \int_y^1{|S(x)|(1-x^2)^{-1/2} \, dx} \cr
\leq & \|S\|_{[0,1]} \int_y^1{(1-x^2)^{-1/2}\, dx} \cr 
\leq & \|S\|_{[0,1]} \int_{\arccos y}^{\arccos 1}{\frac{-\cos t}{\cos t} \, dt} 
= \|S\|_{[0,1]} \int_{0}^{\arccos y}{\, dt} \cr 
\leq & (2(1-y))^{1/2}\|S\|_{[0,1]} \leq \left( \frac{20k}{n} \right)^{1/2} \|S\|_{[0,1]}\,. \cr \endsplit \tag 5.6$$
In the third line of (5.6) we used the substitution $x = \cos t$, while in the fourth line (5.6) we 
used the inequality $\arccos y \leq (2(1-y))^{1/2}$ which follows from the inequality 
$\cos \tau \geq 1- \tau^2/2$ with $\tau = \arccos y$. (Note also that (5.5) and (5.6) show 
that in the sum on the right-hand side of (5.4) the second term is the dominating one.)  
Combining (5.4), (5.5), and (5.6) we obtain
$$V_0^1(P) = \int_0^1{|P^{\prime}(x)| \, dx} \leq (kn)^{-1/2}\|S\|_{[0,1]} + \left( \frac{20k}{n} \right)^{1/2}\|S\|_{[0,1]}  
\leq 6(k/n)^{1/2}\|S\|_{[0,1]}\,,$$ 
and the lower bound of Theorem 2.2 follows.
\qed \enddemo

\demo{Proof of the upper bound in Theorem 2.1}
If $1 \leq k \leq 5$, then the upper bound of the theorem follows by considering $P \in {\Cal P}_{n,k}$ 
defined by $P(x) = x^{n+1}$. So we can assume that $k \geq 6$. Without loss of generality we may assume that 
$n = 2\nu \geq 0$, $k = 2\kappa + 2 \geq 6$ are even, and $20\kappa \leq \nu$.
Let $T := T_{\nu,\kappa}$ be the Chebyshev polynomial for the M\"untz space
$$\text{\rm span}\{x^{\nu}, x^{\nu+1}, \ldots, x^{\nu+\kappa}\}$$
on $[0,1]$ normalized so that $T(1)=1$.
We define $P \in {\Cal P}_{n+k-1}$ of the form 
$$P(x) = x^{n+1}Q(x)\,, \qquad Q \in {\Cal P}_{k-2}\,,$$ 
by 
$$P(x) = \int_0^x{T(u)^2 \, du}\,.$$
Lemma 3.9 implies that 
$$P(1) = |P(1)| \geq \frac{c_5\kappa}{\nu}\,. \tag 5.7$$
Observe that 
$$|P^{\prime}(y)| = T(y)^2 \leq 1\,, \qquad y \in [0,1]\,,$$
and hence  
$$\|P^{\prime}\|_{[0,1]} \leq 1\,. \tag 5.8$$
Combining (5.7) and (5.8) we have
$$\frac{\|P^{\prime}\|_{[0,1]}}{|P(1)|} \leq \frac{1}{c_5\kappa/\nu} = \frac{1}{c_5}\frac{\nu}{\kappa} \leq \frac{c_8n}{k}$$  
with an absolute constant $c_8 > 0$, and the upper bound of Theorem 2.1 follows. 
\qed \enddemo

\demo{Proof of the upper bound in Theorem 2.2}
If $1 \leq k \leq 5$, then the upper bound of the theorem follows by considering $P \in {\Cal P}_{n,k}$
defined by $P(x) = x^{n+1}$. So we can assume that $k \geq 6$. Without loss of generality we may assume that
$n = 2\nu \geq 0$, $k = 2\kappa + 2 \geq 6$ are even, and $20\kappa \leq \nu$.
We define $P \in {\Cal P}_{n+k-1}$ of the form 
$$P(z) = x^{n+1}Q(x)\,, \qquad Q \in {\Cal P}_{k-2}\,,$$ 
by 
$$P(x) = \int_0^x{T(u)^2 \, du}\,,$$
as in the proof of the upper bound in Theorem 2.1.
Let $y$ be a number such that
$$|P^{\prime}(y)\sqrt{1-y^2}| = \|P^{\prime}(x)\sqrt{1-x^2}\|_{[0,1]}\,.$$
Lemma 3.5 implies that $y \geq 1 - 10k/n$, and hence
$$\|P^{\prime}(x)\sqrt{1-x^2}|\|_{[0,1]} = |P^{\prime}(y)\sqrt{1-y^2}| = T(y)^2 \sqrt{1-y^2} 
\leq \sqrt{1-y^2} \leq \left( \frac{20k}{n} \right)^{1/2}\,. \tag 5.9$$
Combining (5.7) and (5.9) we have
$$\frac{\|P^{\prime}(x)\sqrt{1-x^2}\|_{[0,1]}}{|P(1)|} \leq \frac{(20k/n)^{1/2}}{c_5k/n} 
\leq c_9 \left(\frac nk\right)^{1/2}$$
with an absolute constant $c_9 > 0$, and the upper bound of Theorem 2.2 follows. 
\qed \enddemo

\head 6. Acknowledgement \endhead
The author thanks Szil\'ard R\'ev\'esz for checking the details of the proofs
in this paper and for his suggestions to make the paper more readable.

\enddocument